\documentclass[11pt]{amsart}

\usepackage{latexsym,amssymb,amsfonts}

\begin{document}
\newtheorem{thm}{Theorem}[section]
\newtheorem{prop}[thm]{Proposition}
\newtheorem{lem}[thm]{Lemma}
\newtheorem{cor}[thm]{Corollary}
\newtheorem{ex}[thm]{Example}
\newtheorem{rmk}[thm]{Remark}

\theoremstyle{definition}
\newtheorem{defn}[thm]{Definition}
\newtheorem{remark}[thm]{Remark}
\newtheorem{example}[thm]{Example}
\newtheorem{reduction}[thm]{Reduction}

\newcommand{\Span}{\operatorname{Span}}
\newcommand{\str}{\operatorname{str}}
\newcommand{\Fr}{\operatorname{Fr}}
\newcommand{\fl}{\operatorname{fl}}
\newcommand{\Id}{\operatorname{Id}}
\newcommand{\id}{\operatorname{id}}
\newcommand{\pr}{\operatorname{pr}}
\newcommand{\Supp}{\operatorname{Supp}}
\newcommand{\Hom}{\operatorname{Hom}}

\renewcommand{\Re}{\operatorname{Re}}
\renewcommand{\Im}{\operatorname{Im}}

\newcommand{\bbA}{{\mathbb A}}
\newcommand{\bbC}{{\mathbb C}}
\newcommand{\bbZ}{{\mathbb Z}}
\newcommand{\bbR}{{\mathbb R}}
\newcommand{\bbQ}{{\mathbb Q}}
\newcommand{\bbF}{{\mathbb F}}
\newcommand{\bbN}{{\mathbb N}}

\newcommand{\kbar}{\overline k}

\newcommand{\brokrarr}{\vphantom{\to}\mathrel{\smash{{-}{\rightarrow}}}}
\newcommand{\mo}{\mathopen<}
\newcommand{\mc}{\mathclose>}

\let\to=\longrightarrow

\let\tilde=\widetilde
\let\hat=\widehat
\numberwithin{equation}{section}

\newcommand{\Ker}{\operatorname{Ker}}
\newcommand{\Aut}{\operatorname{Aut}}
\newcommand{\Br}{\operatorname{Br}}
\newcommand{\sdp}{\mathbin{{>}\!{\triangleleft}}}
\newcommand{\Alt}{\mbox{A}}
\newcommand{\GL}{\mbox{\rm GL}}
\newcommand{\SL}{\mbox{\rm SL}}
\newcommand{\rank}{\operatorname{rank}}
\newcommand{\aut}{\operatorname{Aut}}
\newcommand{\Char}{\mbox{\rm char\,}} %% \char is already a command
\newcommand{\diag}{\mbox{\rm diag}}
\newcommand{\Gal}{\operatorname{Gal}}
\newcommand{\galois}{\Gal}
\newcommand{\cd}{\operatorname{cd}}
\newcommand{\lra}{\longrightarrow}
\newcommand{\Mat}{\mbox{\rm M}}
\newcommand{\Mn}{\mbox{\rm M}_n}
\newcommand{\ord}{\mathop{\rm ord}\nolimits}
\newcommand{\Sym}{{\operatorname{S}}}
\newcommand{\tr}{\mbox{\rm tr}}
\newcommand{\trace}{\tr}
\newcommand{\exponent}{\operatorname{exp}}

\tolerance=9999 \hbadness=9999

\title[Trace forms]{Trace forms of Galois extensions
in the presence of a fourth root of unity}
\author[J. MIN\'{A}\v{C}]{J. Min\'{a}\v{c} $^{*}$}
\address{Department of Mathematics, University of Western Ontario,
London, Ontario N6A 5B7, Canada}
\thanks{$^*${Partially supported by an NSERC research grant}}
\email{minac@uwo.ca}
\author[Z. REICHSTEIN]{Z. Reichstein $^{*}$}
\address{Department of Mathematics, University of British Columbia,
Vancouver, BC, Canada V6T 1Z2}
\email{reichst@math.ubc.ca}
%%%%%%%%%%%%%%%%%%%%
% 11E81 Algebraic theory of quadratic forms; Witt groups and rings.
% 12F05 Algebraic field extensions
% 12J10 Valued fields
% 20D15 Nilpotent groups, $p$-groups
% 12J10 Valued fields
%%%%%%%%%%%%%%%%%%%%%

\begin{abstract}
We study quadratic forms that can occur as trace forms
$q_{L/K}$ of Galois field extensions $L/K$, under the assumption
that $K$ contains a primitive $4$th root of unity.
M. Epkenhans conjectured that $q_{L/K}$ is always a scaled
Pfister form. We prove this conjecture and classify
the finite groups $G$ which admit a $G$-Galois
extension $L/K$ with a non-hyperbolic trace form.
We also give several applications of these results.
\end{abstract}

\maketitle

\section{Introduction}
\label{sect.intro}

The trace form of a finite separable field extension (or, more generally
of an \'etale algebra) $L/K$ is the non-degenerate quadratic form
$q_{L/K} \colon x \mapsto \tr_{L/K}(x^2)$ defined over $K$.
In this paper we shall address the following problem: Given
a finite group $G$, which quadratic forms over $K$ are trace
forms of $G$-Galois extensions $L/K$? This question has been
extensively studied; see, e.g.~\cite{dek}
and the references there.  In~\cite{kr} D.-S. Kang
and the second author obtained the following partial answer:

\begin{thm} \label{thm0}
Let $L/K$ be a $G$-Galois extension and let
$S$ be a Sylow 2-subgroup of $G$. Assume

\smallskip
(a) $S$ is not abelian, and

\smallskip
(b) $K$ contains a primitive $e$th root of unity, where
\[ e = \min \{ \exp(H) \, | \, \text{$H$ is
a non-abelian subgroup of $S$} \} \, . \]

\smallskip
\noindent
Then the trace form $q_{L/K}$ is hyperbolic over $K$.
\end{thm}

In this paper we will study trace forms of $G$-Galois extensions $L/K$,
assuming only that $K$ contains a primitive $4$th root of unity.
M. Epkenhans has conjectured that in this situation $q_{L/K}$
is always a scaled Pfister form.  Our first main result
is a proof of this conjecture.  Before giving
the precise statement, we introduce some notations.

If $G$ is a group and $i \geq 1$ is an integer,
we set $G^i = \mo g^i \, | \, g \in G \mc
\, \triangleleft \; G$.
If $S$ is a finite $2$-group, then $S^2 = \Fr(S)$ is the Frattini
subgroup of $S$.
The Frattini rank $r$ of $S$ is the rank of the
elementary abelian group $S/S^2 \simeq
(\bbZ/2)^r$. Note that the Frattini rank of $S$
equals the cardinality of any minimal generating
set of $S$; see, e.g.,~\cite[7.3]{scott}.

\begin{thm} \label{thm1} Suppose $K$ is a field containing a primitive
$4$th root of unity, $L/K$ is $G$-Galois extension, $S$ is
a Sylow 2-subgroup of $G$, and $r$ is the Frattini rank of $S$. Then
the trace form $q_{L/K}$ is Witt-equivalent to the scaled Pfister
form $\mo |S| \mc \otimes \ll a_1, \dots, a_r \gg$,
for some $a_1, \dots, a_r \in K^*$.
\end{thm}

Several remarks are in order, regarding Theorem~\ref{thm1}.
First of all, both Theorem~\ref{thm0} and~\ref{thm1} remain true
for Galois $K$-algebras $L$ that are not necessarily fields.
The reason is that both are enough to check for a single ``versal"
$G$-Galois algebra, which is a field;
cf. e.g.,~\cite[Proposition 2.5]{kr}.

Secondly, Theorem~\ref{thm1} was previously known
for $|S| \leq 16$; see~\cite[Corollary 6, p. 227]{dek}.

Thirdly, the ``scaling factor" of $\mo |S| \mc$ presents
only a minor inconvenience in working with the trace
form $q_{L/K}$. It can be dropped if $|S|$ is a square
in $K$ (and, in particular, if $K$ contains a primitive
$8$th root of unity; cf. Remark~\ref{rem1.modular})
and replaced by $\mo 2 \mc$ in all other cases.

Finally, the requirement that $K$ should contain a primitive
$4$th root of unity is essential. Indeed, let
$K = \bbQ$ and $L = \bbQ(\sqrt{2 + \sqrt{2}})$.
By~\cite[Proposition 8]{dek} (with $q = a = b = 1$ and $D = 2$),
the field extension $L/K$ is Galois, with
$\Gal(L/K) = \bbZ/4$ and the trace form
$q_{L/K} = \mo 1, 2, 1, 1\mc$. This form is positive-definite and thus
anisotropic. Consequently, $q_{L/K}$ cannot be
Witt-equivalent to a 2-dimensional form. This shows that
Theorem~\ref{thm1} fails for this extension.

\smallskip
Our second main result is a complete description of those
finite groups $G$ which admit a $G$-Galois extension
$L/K$ with a non-hyperbolic trace form. (Here we assume that
$K$ contains a primitive root of unity of degree $2^m$ for a fixed
$m \geq 2$.) It turns out that these groups belong to a rather small
but interesting family that was previously studied
for entirely different reasons.

\begin{thm} \label{thm2} Let $G$ be a finite group,
$S$ be a Sylow 2-subgroup of $G$ and $m \geq 2$ be an integer.
Then the following conditions are equivalent:

\smallskip
(a) there exists a $S$-Galois extension $E/F$ such that $F$
contains a primitive root of unity of degree $2^m$
and the trace form $q_{E/F}$ is not hyperbolic,

\smallskip
(b) there exists a $G$-Galois extension $L/K$ such that
$K$ contains a primitive root of unity of degree $2^m$
and the trace form $q_{L/K}$ is not hyperbolic,

\smallskip
(c) $T/T^{2^m}$ is abelian for every subgroup $T$ of $S$,

\smallskip
(d) there exist an integer
$s \geq m$, an abelian subgroup $A \triangleleft S$,
and an element $t \in S$  such that $S = \mo A, t \mc$ and
$t a t^{-1} = a^{1 + 2^s}$ for every $a \in A$.
\end{thm}

A simple argument based on Sylow's theorem shows that
condition (c) is equivalent to $H/H^{2^m}$ being abelian for every
subgroup $H$ of $G$ (see Remark~\ref{rem.sylow}). Note also that
the $G$-Galois extension $L/K$ in part (b) can be chosen so that
$\Char(K) = 0$ (see Remark~\ref{rem.char0}) and $K$
does not contain a primitive root of unity of degree $2^{m+1}$
(see Remark~\ref{rem.root}).

The $2$-groups $T$ appearing in condition (c)
are {\em powerful} in the sense of Lubotzky
and Mann~\cite{lm}.  Their results on
the structure of powerful groups will be used in the proof 
of Theorem~\ref{thm2}, along with theorems 
of Iwasawa~\cite{iwasawa} and Engler-Koenigsmann~\cite{ek}.

Theorems~\ref{thm1} and~\ref{thm2} have a natural cohomological
interpretation. 
Let $G$ be a finite group, $S$ be a Sylow $2$-subgroup
of $G$, $r$ be the Frattini rank of $S$ and $K$ be a field containing 
a primitive root of unity of degree $2^m$ for some integer $m \geq 2$. 
Then to every $G$-Galois field extension $L/K$ (and, more generally, to
a $G$-Galois $K$-algebra $L$) we can associate the well-defined
cohomology class $\phi(L) = (a_1) \cdot (a_2) \dots (a_r)$ in
$H^r(K, \bbZ/2 \bbZ)$, where $a_1, \dots, a_r$ are as in Theorem~\ref{thm1}.
In other words, $\phi(L)$ is the Arason invariant of the Pfister form
$\mo |S| \mc \otimes q_{L/K}$; cf.~\cite[Section 1]{arason}.
The map $\phi$ so defined is 
easily seen to be a cohomological invariant
\[ \phi \colon H^1(\ast, G) \lra H^r(\ast, \bbZ/2 \bbZ) \, ,\]
where $\ast$ ranges over the category of fields containing a primitive 
$2^m$th root of unity.
(Recall that the non-abelian cohomology set $H^1(K, G)$
parametrizes $G$-Galois algebras over $K$.) Theorem~\ref{thm2} 
gives equivalent conditions for this cohomological invariant to be 
non-trivial.

The rest of this paper is structured as follows.
Theorem~\ref{thm1} is proved in Sections~\ref{sect.serre}
and~\ref{sect.thm1}.  Theorem ~\ref{thm2} is proved
in Sections~\ref{sect.iwasawa} - \ref{sect.pf3}.
In Section~\ref{sect.appl} we discuss a number
of applications of these results.  In particular,
we show that the trace form of a $G$-Galois field
extension $L/K$  is hyperbolic if the field $K$
is ``sufficiently small" in a suitable sense
(see Proposition~\ref{prop2}) or if $G$ is a simple
group whose Sylow $2$-subgroups are non-abelian
(see Proposition~\ref{prop7.8}). In the last section
we give a description of quadratic forms that can occur
as trace forms of $M(2^n)$-Galois extensions, where
$$M(2^{n})=<\sigma,\tau|\sigma^{2^{n-1}}=1=\tau^{2},
\tau\sigma\tau=\sigma^{1 + 2^{n-2}}>.$$

\section*{Acknowledgments}
We would like to thank the organizers of the MSRI programs
in Galois Theory and Noncommutative Algebra and MSRI staff
for giving us an opportunity to meet and to begin our collaboration
in the Fall of 1999. We are also grateful to M.~Epkenhans, 
A.~Pfister and J.-P.~Serre for stimulating
discussions, and to G.~Bhandari, J.-P.~Serre 
and the anonymous referee for helpful comments 
on an earlier draft of this paper.

\section{Orthogonal 2-groups}
\label{sect.serre}

Most of our subsequent results will be based on the following lemma,
communicated to us by J.-P. Serre.

\begin{lem} \label{lem.serre}
Let $K$ be a field containing a primitive $4$th root of unity,
$(V, q)$ be a non-degenerate finite-dimensional
quadratic space over $K$ and $G$ be a finite $2$-subgroup,
acting orthogonally on $V$.  Then $V$ can be decomposed as an
orthogonal sum $V = V^{\Fr(G)} \oplus V_0$, such that the
restriction of $q$ to $V_0$ is hyperbolic.
\end{lem}

Here, as usual, $V^{\Fr(G)} = \{ v \in V \, | \, h(v) = v$ for every 
$h \in \Fr(G) \}$, and we allow the trivial hyperbolic 
quadratic space $V_0 = \{0\}$.

\begin{proof} We argue by induction on $\dim(V) + |G|$.
Assume, to the contrary, that the lemma fails for some $V$, $q$ and $G$;
choose a counterexample with $\dim(V) + |G|$ as small as possible. Then
$G$ acts faithfully on $V$; otherwise we could obtain a counterexample
with a smaller value of $\dim(V) + |G|$ by keeping the same $V$ and
replacing $G$ by $G/N$, where $N$ is the kernel of this action.

We claim that every index 2 subgroup of $G$ is elementary abelian.
Indeed, assume the contrary: $\Fr(H) \neq \{ 1 \}$ for some index 2
subgroup $H$.  Equivalently,  $V^{\Fr(H)} \neq V$.  
Since $| H | + \dim V < |G| + \dim V$,
our induction hypothesis applies and we can write $V$ as an orthogonal sum
\[ V=V^{\Fr(H)}\oplus V_1 \, , \]
where the restriction of $q$ to $V_1$ is hyperbolic. 
In particular, $(V^{\Fr(H)}, \, q_{| \, V^{\Fr(H)}})$ 
is a regular quadratic space;
see~\cite[p. 11, Corollary 2.6]{lam}.  Since $\Fr(H)$ is a normal 
subgroup of $G$, the action of $G$ restricts to $V^{\Fr(H)}$.
This restricted action is once again orthogonal, and since 
$\dim V^{\Fr(H)}<\dim V$, we can apply our induction assumption
to write $V^{\Fr(H)}$ as an orthogonal sum
$$V^{\Fr(H)}=V^{\Fr(G)}\oplus V_2 \, ,$$
where the restriction of $q$ to $V_2$ is hyperbolic. To sum up,
\[ V=V^{\Fr(H)}\oplus V_1 = V^{\Fr(G)}\oplus V_0 \, , \]
where the restriction of $q$ to $V_0 = V_1 \oplus V_2$ is hyperbolic, 
contradicting our choice of $V$ and $G$. This contradiction proves the claim.

If every element of $G$ has order $\leq 2$ then $G$ is itself 
elementary abelian. In this case the lemma is trivial,
because $\Fr(G) = \{ 1 \}$. Thus we may assume $G$ has an element $g$
of order $4$. By the claim we just proved, $g$ is not contained in any
subgroup of $G$ of index $2$. In other words, $\mo g \mc$ is 
not contained in any proper 
subgroup of $G$, i.e., $G = \mo g \mc \simeq \bbZ/4$. 
We shall thus concentrate on this case for the rest of the proof. 
Note that $\Fr(G) = \mo g^2 \mc$. We now proceed with an explicit
description of $V_0$.

Now recall that $K$ is assumed to contain a primitive $4$th root of unity;
we will denote it by $\zeta$.  Since $g^4 = 1$, we can decompose 
$V$ as a direct sum of the four eigenspaces for $g$:
\begin{equation} \label{eigenspaces}
V = V_1 \oplus V_{-1} \oplus V_{\zeta} \oplus V_{-\zeta} \, ,
\end{equation}
where $V_{\alpha} = \{ v \in V \, | \, g(v) = \alpha v \}$.
Note that if $x \in V_{\alpha}$ and $y \in V_{\beta}$ then
\[ B(x,y)=B(g(x), g(y)) = \alpha \beta B(x, y) \]
and thus
\begin{equation} \label{e.hyperb}
\text{$B(x, y) = 0$ whenever $\alpha \beta \neq 1$.}
\end{equation}
Here $B$ denotes the bilinear form associated with the quadratic form $q$.

In particular $V^{\Fr(G)} = V_1 \oplus V_{-1}$ is orthogonal to
$V_{\zeta} \oplus V_{-\zeta}$, and thus we can take
$V_0 = V_{\zeta} \oplus V_{-\zeta}$. By~\eqref{e.hyperb} both
$V_{\zeta}$ and $V_{-\zeta}$ are totally
isotropic. Thus $V_0$ contains a totally isotropic space of dimension
at least half the dimension of $V_0$. Observe also that from~\eqref{e.hyperb},
and from our assumption that $q$ is non-degenerate on $V$, it follows that
$q$ is non-degenerate on $V_0$. Thus we see that $V_0$ is hyperbolic;
see \cite[Chapter~1,~Theorem 3.4(i)]{lam}. To sum up,
\[ V = (V_1 \oplus V_{-1}) \oplus (V_{\zeta} \oplus V_{-\zeta}) =
 V^{\Fr(G)} \oplus V_0 \, , \]
where the restriction of $q$ to $V_0$ is hyperbolic. This contradicts 
our choice of $G$ and $V$, thus completing the proof of Lemma~\ref{lem.serre}.
\end{proof}

\begin{cor} \label{cor.serre} Let $G$ be a finite $2$-group and
$L/K$ be a $G$-Galois extension. Assume $K$ contains a primitive
$4$th root of unity. Then

\smallskip
(a) $q_{L/K} \simeq \mo |\Fr(G)| \mc \otimes q_{L^{\Fr(G)}/K}$.

\smallskip
(b) More generally, for any normal subgroup $H \subset \Fr(G)$,
\[ q_{L^H/K} \simeq \mo [\Fr(G): H] \mc \otimes q_{L^{\Fr(G)}/K} \, . \]

\smallskip
\noindent
Here $\simeq$ denotes Witt equivalence.
\end{cor}

\begin{proof}
(a) The 2-group $G$ acts orthogonally on the
quadratic space $(V = L, q_{L/K})$ over $K$. By
Lemma~\ref{lem.serre}, $q_{L/K}$ is Witt-equivalent
to its restriction to $L^{\Fr(G)}$. Finally, for every
$x \in L^{\Fr(G)}$, we have 
\[ q_{L/K}(x) = |\Fr(G)| \,  q_{L^{\Fr(G)}}(x) \, , \]
and part (a) follows.

(b) Apply part (a) to the $G/H$-Galois extension
$L^H/K$, remembering that $\Fr(G/H) = \Fr(G)/H$.
%  The 2-group $G$ acts orthogonally on the
% quadratic space $(L^H, q_{L^H/K})$ over $K$. By
% Lemma~\ref{lem.serre}, $q_{L/K}$ is Witt-equivalent to the form $q_0$
% on $L^{\Fr(G)}$ given by $q_0(x) = q_{L^H/K}(x)$ for every $x \in L^{\Fr(G)}$.
% Then \[ q_0(x) = [\Fr(G): H] \, \tr_{L^{\Fr(G)}/K}(x^2) \, , \] and part (a)
% follows. To deduce (b) from (a), set $H = \{ 1 \}$.
\end{proof}

\section{Conclusion of the proof of Theorem~\ref{thm1}}
\label{sect.thm1}

As usual, given $a_1,a_2,\dots,a_n\in K^{*}$,
$\ll a_1,\dots,a_n \gg = \otimes_{i=1}^{n}\mo 1,-a_i \mc$
will denote an $n$-fold Pfister form. Note that since we always
assume $K$ contains a primitive $4$th root of unity,
$$\ll a_1,\dots,a_n \gg \simeq \otimes_{i=1}^{n}\mo 1, a_i \mc \, . $$

We now begin the proof of Theorem~\ref{thm1}
by reducing to the case where $G = S$ is a 2-group.

\begin{lem} \label{lem2.1} Let $G$ be a finite group, $K$ be
a field containing a primitive $4$th root of unity, $L/K$
be a $G$-Galois extension, $S$ be the Sylow 2-subgroup of $G$,
$K_1 = L^S$ and $\phi \colon W(K) \lra W(K_1)$
be the natural (extension of scalars) homomorphism
of Witt rings.

\smallskip
(a) $($cf. \cite[6.1.1]{bs}$)$
$q_{L/K_1} = \phi(q_{L/K})$ in $W(K_1)$.

\smallskip
(b) $q_{L/K}$ is hyperbolic if and only if $q_{L/K_1}$ is hyperbolic.

\smallskip
(c) Let $a \in K^*$. Then $q_{L/K} =
\mo a \mc \otimes \ll a_1, \dots, a_r \gg$ in $W(K)$,
for some $a_1, \dots, a_r \in K^*$, if and only if
$q_{L/K_1} = \mo a \mc \otimes \ll b_1, \dots, b_r \gg$ in $W(K_1)$
for some $b_1, \dots, b_r \in K_1^*$.
\end{lem}

\begin{proof} (a) $\phi(q_{L/K})$ is clearly the trace form
of the $K_1$-algebra
$L_1 = L \otimes_{K} K_1$ and $L_1$ is isomorphic, as a
$K_1$-algebra, to
\begin{equation} \label{e2.1}
\text{$L \oplus \dots \oplus L$ ($m$ times),}
\end{equation}
where $m = [G: S]$ is odd. Moreover,~\eqref{e2.1} is an orthogonal
direct sum with respect to the trace form. Thus
\[ \text{$\phi(q_{L/K}) = q_{L/K_1} \oplus \dots \oplus q_{L/K_1}$
($m$ times);} \]
cf.~\cite[Theorem I.5.1]{cp}.
Since we are assuming $K$ (and thus $K_1$) contains a primitive
4th root of unity, $2W(K_1)=\{0\}$, and part (a) follows.

\smallskip
By Springer's theorem, $\phi$ is injective;
see, e.g.,~\cite[Theorem 7.2.3]{lam}. Part
(b) now follows from (a).

\smallskip
(c) By Rost's theorem on the descent of Pfister
forms~\cite[Section 3]{rost} (see also~\cite[4.4.1]{bs}),
$\mo a \mc \otimes q_{L/K}$ is Witt-equivalent to
a Pfister form over $K$ if and only if
$\mo a \mc \otimes q_{L/K_1}$ is Witt-equivalent to
a Pfister form over $K_1$.
\end{proof}

We now continue with the proof of Theorem~\ref{thm1}.
By Lemma~\ref{lem2.1}(c) we may assume that $G$ is a 2-group.
By Corollary~\ref{cor.serre}
\[ q_{L/K} \simeq \mo |\Fr(G)| \mc \otimes q_{L^{\Fr(G)}/K} \, . \]
Note that $L^{\Fr(G)}/K$ is a $G/\Fr(G)$-Galois extension, where
$G/\Fr(G) \simeq (\bbZ/2)^r$. Thus it is enough to prove
Theorem~\ref{thm1} in the case where $\Gal(L/K)$ is an elementary
abelian 2-group; indeed, if we know that
\[ q_{L^{\Fr(G)}/K} \simeq \mo |G/\Fr(G)| \mc \otimes
\text{($r$-fold Pfister form)} \, . \] then 
by Corollary~\ref{cor.serre}(a) 
\[ q_{L/K} \simeq \mo | \Fr(G) | \mc \otimes q_{L^{\Fr(G)}/K} \simeq
\mo |G| \mc \otimes \text{($r$-fold Pfister form),} \]
as claimed.

Now assume $G = (\bbZ/2)^r$.
Here any $G$-Galois extension $L/K$ has the form
$L = K(\sqrt{a_1}, \dots, \sqrt{a_r})$, for some
$a_1, \dots, a_r \in K^*$, and an easy computation in the basis
$\{ a_1^{\frac{\epsilon_1}{2}}\dots a_r^{\frac{\epsilon_r}{2}} \}$,
with $\epsilon_1, \dots, \epsilon_r = 0, 1$, shows that
\begin{equation} \label{e.elem-abelian}
q_{L/K} \simeq \mo 2^r \mc \otimes \ll a_1, \dots, a_r \gg \, ;
\end{equation}
cf.~\cite[6.2.1]{bs} or~\cite[Lemmas 2.1(b) and 2.2]{kr}.
This completes the proof of Theorem~\ref{thm1}.
\qed

\section{Iwasawa structures}
\label{sect.iwasawa}

An {\em Iwasawa structure of level} $s \geq 1$ on a $2$-group $G$
is a normal abelian subgroup $A$
and an element $t$ such that $G = \mo A, t \mc$ and
\[ \text{$t a t^{-1} = a^{1 + 2^s}$ for every $a \in A$.} \]
Informally speaking, the higher the level is,
the closer $G$ is to an abelian group.
In particular, if $\exp(A) = 2^e$ and $s \geq e$
then $G$ is abelian. Conversely, any finite abelian
$2$-group $G$ of exponent $\leq 2^s$ admits
an Iwasawa structure of level $s$,
with $A = G$ and $t = \{ 1 \}$.

If a $2$-group $G$ admits an Iwasawa structure of level $\ge 2$, we will
call $G$ an Iwasawa group. Note that the level
of an Iwasawa group $G$ is not well-defined
in general, since $G$ may admit Iwasawa
structures of different levels (see
Example~\ref{ex.iwasawa} below).

For any $2$-group $G$ we define the {\em strength} of $G$ by
\[ \str(G) = \max \;  \{ m \; | \; \text{$G/G^{2^m}$
is abelian} \} \, . \]
In particular, $\str(G) = \infty$ iff $G$ is abelian and
$\str(G) \geq 2$ iff $G$ is {\em powerful}
in the sense of Lubotzky and Mann; cf.~\cite[Definition, p. 499]{lm}.

\begin{lem} \label{lem4.3} Suppose that $G$ is
a finite $2$-group which admits an Iwasawa structure
$(A,t)$ of level $s$. Then

\smallskip
(a) $[G, G] = A^{2^s}$,

\smallskip
(b) $\str(G) \geq s$,

\smallskip
(c) If $s \geq 2$ then $G^{2^m} = \mo A^{2^m}, t^{2^m} \mc$ for every $m \in
\bbN$.
\end{lem}

\begin{proof} (a) From the definition of an Iwasawa structure of level $s$,
we see that $A^{2^s} \subset [G, G]$ and $G/A^{2^s}$ is abelian.
Hence, $[G, G] = A^{2^s}$.

\smallskip
(b) By part (a) $G/A^{2^s}$ is commutative. Hence, so is
$G/G^{2^s}$, and part (b) follows.

\smallskip
(c) By part (b), $\str(G) \geq 2$. Thus
$[G, G] \subset G^4$, i.e., $G$ is a powerful $2$-group.
The desired conclusion now follows from~\cite[Theorem~2.7]{ddms}.
\end{proof}

We remark that part (c) remains true even if $s = 1$. This stronger assertion
will not be used in the sequel; we leave it as an exercise for the reader.

\begin{example} \label{ex.iwasawa}
The inequality $\str(G) \geq s$ may be strict, even if $G$
is non-abelian.  Indeed, let
$$G=\mo a,t\vert a^{2^5}=1,a^{2^2}=t^{2^3},tat^{-1}=a^{1+8}\mc.$$
One checks readily that $G$ is a metacyclic group of
order $2^8$ and that $G$ admits an Iwasawa structure $(A,t)$
of level $3$, where $A=\mo a \mc$. We claim that
$\str(G) = 4$. By Lemma~\ref{lem4.3},
$[G,G]=\mo a^8 \mc$. Since
$a^8=t^{16}$, we see that $[G,G]$ is contained in $G^{16}$ but not
in $G^{32} = \mo a^{32}, t^{32} \mc = \mo a^{16} \mc$. Thus
$\str(G) = 4$, as claimed.

On the other hand, observe that $G$ admits another Iwasawa structure
$(\tilde{A},\tilde{t})$ of level $4$, where $\tilde{A}=\mo t \mc$ and
$\tilde{t}=a^{-1}$. Indeed,
have $\tilde{t} \, t \, \tilde{t}^{-1}=a^{-1}ta=t^{1+2^4}$.
Thus we see that by switching the role of $t$ and $a^{-1}$,
we are able to find another Iwasawa structure whose level
equals the strength of $G$. In the next lemma we shall
show that such a switch is always possible.
\end{example}

\begin{lem} \label{lem4.5} Suppose $G$ be a non-abelian
Iwasawa $2$-group.  Then
\[ \str(G)= \max \{ \text{level}(A,t) \} \, , \]
where the maximum is taken over all Iwasawa
structures $(A,t)$ on $G$.
\end{lem}

\begin{proof} Let $m = \str(G)$ and
$(A,t)$ is an Iwasawa structure on $G$ of level $s$. By
Lemma~\ref{lem4.3}, $s \leq m$. If $s = m$ we are done.
Thus we may assume $s < m$. Our goal is to construct another
Iwasawa structure on $G$ of level $m$.

Since $G$ is an Iwasawa $2$-group, $m \geq 2$. Thus
$[G, G] \subset G^4$, so that $G$ is a powerful group.
By Lemma~\ref{lem4.3},
\[ A^{2^s} = [G, G] \subset G^{2^m} = \mo A^{2^m},t^{2^m}\mc \, . \]
We now see that the group $G^{2^m}/A^{2^m}$ is cyclic, and hence,
so is its subgroup $A^{2^s}/A^{2^m}$. Since $s < m$ this implies
that $A^{2^s}$ is itself cyclic.

Let $a^{2^s}=t^{2^m}$ be a generator of $A^{2^s}$ with
$a\in A$. Since the order of $a$ is equal to the exponent
of $A$, we see that there exists a subgroup $B$ of $A$ such that
$A=\mo a\mc\oplus B$. Moreover, since $A^{2^s}/A^{2^m}$ is
cyclic, we see that $B^{2^s}=\{1\}$. Therefore,
$tbt^{-1}=b^{1+2^s}=b$ for each $b\in B$, and
$B$ is a subgroup of the center $Z(G)$ of $G$.

Set $\tilde{A}=\mo t,B \mc$ and $\tilde{t}=a^{-1}$. We claim that
$(\tilde{A},\tilde{t})$ is an Iwasawa structure on $G$ of level $m$.
First we have
$$\mo\tilde{A},\tilde{t}\mc=\mo t,B,a^{-1}\mc=\mo t,A\mc=G.$$
Also $\tilde{A}$ is an abelian subgroup of $G$ as $B\subset Z(G)$.
Further $\tilde{t} \, t \, \tilde{t} \, ^{-1} =
a^{-1}ta=a^{-1}a^{1+2^s}t=a^{2^s}t=t^{1 + 2^m}$,
as $a^{2^s}=t^{2^m}$. Because $\tilde{A}=\mo B,t\mc$ and $B\subset Z(G)$
we see that $\tilde{t} \, \tilde{a} \,  \tilde{t} \, ^{-1} =
\tilde{a} \, ^{1+2^m}$
for each $\tilde{a}\in\tilde{A}$. Hence $(\tilde{A},\tilde{t})$
is the Iwasawa structure of level $m$.
\end{proof}

\begin{remark} \label{rem.iwasawa}
In view of Lemma~\ref{lem4.5}, a $2$-group $S$ satisfies
condition (d) of Theorem~\ref{thm2} if and only if it is an
Iwasawa group of strength $\geq m$.
\end{remark}

\section{Proof of Theorem~\ref{thm2} (a) $\Longrightarrow $ (b)
$\Longrightarrow$ (c) $\Longrightarrow$ (d)}
\label{sect.pf1}

\smallskip
(a) $\Longrightarrow $ (b):
Let $k$ be the subfield of $F$ generated by the prime field and
the primitive $2^m$th root of unity
and let $V$ be a faithful linear representation of $G$ over
$k$ (e.g., we can take $V$ to be the group algebra $k[G]$).
Denote the field of rational functions on $V$ by $k(V)$.
Since the trace form of the $S$-Galois extension $E/F$ is not
hyperbolic~\cite[Proposition 2.5]{kr} tells us that
the trace form of $k(V)/k(V)^S$ is not hyperbolic. Now by
Lemma~\ref{lem2.1}(b), $k(V)/k(V)^G$ is not
hyperbolic either. Thus we can take $L = k(V)$ and $K = k(V)^G$.

\smallskip
(b) $\Longrightarrow $ (c): Let $L/K$ be a $G$-Galois field
extension with a non-hyperbolic trace form, as in (b). Assume,
to the contrary, that $T/T^{2^m}$ is non-abelian for some
subgroup $T$ of $S$. Then the trace form of $L/L^{T}$
is still non-hyperbolic; see \cite[Lemma 2.1(c)]{kr}. Thus,
replacing $G$ by $T$ and $K$ by $L^T$, we may assume
$G = T$.

Now let $H = G^{2^m}$. Then $L^H/K$ is a Galois extension with Galois group
$G/H$, which by our assumption, is non-abelian of exponent $\leq 2^m$.
Thus, by Theorem~\ref{thm0}, $q_{L^H/K}$ is hyperbolic.
Now, since $H \subset G^2 = \Fr(G)$, Corollary~\ref{cor.serre}
tells us that $q_{L/K}$ is hyperbolic as well, contradicting our assumption.

\smallskip
(c) $\Longrightarrow $ (d): By our assumption every subgroup $T$ of $S$
satisfies $[T, T] \subset T^4$, i.e., $T$ is powerful.
By~\cite[Theorem 4.3.1]{lm} this implies that $S$ is modular
but not Hamiltonian. On the other hand, by a theorem of
Iwasawa~\cite{iwasawa} modular non-Hamiltonian $2$-groups are
precisely the $2$-groups that admit an Iwasawa structure of
of level $\ge 2$.
\footnote{The proofs of Iwasawa's theorem in~\cite{iwasawa} 
and~\cite[Theorem 14]{suzuki} had some gaps that were 
later pointed out and closed by Napolitani~\cite{napolitani}. 
For a detailed exposition of Iwasawa's theorem and related 
group-theoretic results, we refer the reader to~\cite{schmidt}.}

It remains to show that $S$ admits an Iwasawa structure 
of level $s \geq m$.
First suppose $S$ is abelian. Then, as we pointed out
in Section~\ref{sect.iwasawa},
we can take $A = S$, $t = 1$, and $s = \max \{ m, e \}$,
where $e$ is the exponent of $S$. Now assume $S$ is not abelian.
Then by our assumption (c), $\str(S) \geq m$. The desired
conclusion now follows from Lemma~\ref{lem4.5}.
\qed

\begin{remark} \label{rem.root}
{\em If $\Char(K) = 0$ then the $G$-Galois extension $L/K$ in part (b) 
can be chosen so that
$K$ does not contain a primitive root of unity of degree $2^{m+1}$.}

\begin{proof}
Let $k = \bbQ(\zeta_{2^m})$ be the subfield of $K$ generated by its
prime subfield and a primitive $2^m$th root of unity.
Let $V = k^n$ be a faithful $G$-representation (over $k$), as in the proof
of the implication $(a) \Longrightarrow (b)$.
Since the trace form of the $S$-Galois extension $E/F$ is not
hyperbolic,~\cite[Proposition 2.5]{kr} tells us that
the trace form of $k(V)/k(V)^S$ is not hyperbolic. Thus we can replace
$E$ by $E' = k(V)$ and $F$ by $F' = k(V)^G$.  Since $k$
is algebraically closed in $E'$, $E'$ (and hence, $F'$)
does not contain a primitive root of unity of degree $2^{m+1}$.
\end{proof}

The same argument goes through in characteristic $p$, provided that
$k = \bbF_{p}(\zeta_{2^m})$ does not contain $\zeta_{2^{m+1}}$.
\end{remark}

\begin{remark} \label{rem.sylow}
{\em Condition (c) of Theorem~\ref{thm2} is equivalent to

\smallskip
(c$'$) $H/H^{2^m}$ is abelian for every subgroup $H$ of $G$.}

\begin{proof}
Clearly, (c$'$) $\Longrightarrow$ (c). To prove the converse,
let $T$ be a Sylow $2$-subgroup of $H$. After replacing $S$
by a conjugate Sylow subgroup in $G$, we may assume $T \subset S$.
Let $\overline{T}$ be the image of $T$ in $H/H^{2^m}$.
We claim that $\overline{T} = H/H^{2^m}$. Indeed,
on the one hand, the exponent of $H/H^{2^m}$ divides $2^m$, so that
$H/H^{2^m}$ is a $2$-group. On the other hand, since $T$ is 
a Sylow $2$-subgroup of $H$, the index $[H:T]$ is odd. The index 
of $\overline{T}$ in $H/H^{2^m}$ is thus 
both odd and a power of $2$; hence, $\overline{T} = H/H^{2^m}$, 
as claimed.

Consequently,
\[ T/T^{2^m} \stackrel{\text{onto}}{\lra}
T/(T \cap H^{2^m}) \simeq H/H^{2^m} \, . \]
If $T/T^{2^m}$ is abelian, then so is $H/H^{2^m}$. This shows that
(c) $\Longrightarrow$ (c$'$).
\end{proof}
\end{remark}

\begin{remark} \label{rem.lattice} Let $G$ be a finite group.
If $A$ and $B$ are subgroups of $G$, we shall denote the set of
intermediate subgroups $A \subset X \subset B$ by $[A, B]$.
This set is naturally a lattice, where $X \wedge Y = X \cap Y$ and
$X \vee Y$ = subgroup generated by $X$ and $Y$.

\smallskip
{\em Let $S$ be a Sylow $2$-subgroup of $G$.
Suppose for some subgroups $A$ and $B$ of $S$,
the map $\varphi_{A, B}:[A, \, A \vee B]\lra[A \wedge B, \, B]$,
defined by $\varphi_{A, B}(X)= A \wedge X$, is not a lattice isomorphism.
Then the trace form $q_{L/K}$ is hyperbolic for every $G$-Galois extension
$L/K$ such that $K$ contains a primitive $4$th root of unity.}

\smallskip
\begin{proof}
If $\varphi_{A, B}$ is not a lattice isomorphism for some $A$ and $B$
then the lattice $[\{ 1 \}, S]$ is not modular;
see~\cite[Theorem~2.1.5]{schmidt}.  Then, by Iwasawa's theorem
(the easy direction), $S$ does not satisfy condition (d)
of Theorem~\ref{thm2}. The desired conclusion follows from the implication
(b) $\Longrightarrow$ (d).
\end{proof}
\end{remark}

\section{Proof of Theorem~\ref{thm2} (d) $\Longrightarrow $ (a):
Preliminary reductions}
\label{sect.pf2}

We begin by observing that for the purpose of proving the implication 
$(d) \Longrightarrow (a)$, we may assume that $G = S$ 
is a $2$-group and that $m = s$.  We shall say that 
$S$ admits a non-hyperbolic trace form if it satisfies 
condition (a) of Theorem~\ref{thm2}.

It is easy to see that every abelian $2$-group admits a non-hyperbolic
trace form; see, e.g.,~\cite[Remark 3.2]{kr}. Thus we will assume
from now on that $S$ is non-abelian.  
Recall that by our assumption (d), $S = \mo A, t \mc$, where $A$ 
is abelian and
\begin{equation} \label{e.iwasawa}
\text{$tat^{-1} = a^{1 + 2^s}$ for every $a \in A$.}
\end{equation}

Our proof of the implication (d) $\Longrightarrow $ (a) of
Theorem~\ref{thm2} will consist of two parts. In this section we
will reduce the problem to the case where 
$A = (\bbZ/2^{e} \bbZ)^r$ and $S$ is a semidirect product 
of $A$ and $\mo t \mc$; in the next section we will show that every $S$ 
of this form admits a non-hyperbolic trace form. (Note that here
$r$ is the Frattini rank of $A$; the Frattini rank of $S$ is $r + 1$.) 

In order to facilitate working with Iwasawa groups, we will
write them in terms of generators and relations.  Decompose 
the abelian $2$-group
\[ A = \mo a_1 \mc \times \dots \times \mo a_r \mc \simeq
\bbZ/2^{e_1} \bbZ \times \dots \times \bbZ/2^{e_r} \bbZ \, ,\]
as a product of cyclic subgroups, where $a_i$ has order $2^{e_i}$. Then
$\exp(A) = 2^e$, where $e = \max \{ e_1, \dots, e_r \}$.
Since $S$ is non-abelian, 
\begin{equation} \label{e.s} s < e \, .
\end{equation}
Denote the order of the image of $t$ in $G/A$ by $2^q$
and let $a_0 = t^{2^q} \in A$.
Note that the order of $a_0$ in $A$ is $2^{-q} |\mo t \mc|$ and,
since $a_0$ commutes with $t$, $a_0^{2^s} = 1$ in $A$.

\begin{lem} \label{lem.generators}
(a) The group
$X = A * \mo t \mc / \mo t a t^{-1} = a^{1 + 2^s} \, | \, a \in A \mc$
is isomorphic to $A \sdp \mo t \mc$, with the action of $t$ on $A$ given
by~\eqref{e.iwasawa}. Here $A * \mo t \mc$ denotes
the free product of the subgroups $A$ and $\mo t \mc$ of $G$.

\smallskip
(b) Let $c \in A$ be an element of order
$2^{-q} |\mo t \mc|$, satisfying $c^{2^s} = 1$ and
$Y = A * \mo t \mc / \mo t^{2^q} = c \, , \;
t a t^{-1} = a^{1 + 2^s} \, | \, a \in A \mc$.
Then every element of $Y$ can be uniquely written in the form
$a t^i$ for some $a \in A$ and $0 \le i < 2^q$.

\smallskip
(c) $S$ is isomorphic to
$Z = A * \mo t \mc / \mo t^{2^q} = a_0 \, , \;
t a t^{-1} = a^{1 + 2^s} \, | \, a \in A \mc$.
\end{lem}

\begin{proof}
(a) Consider the natural surjective homomorphism $X \lra A \sdp \mo t \mc$,
taking $a$ to $a$ and $t$ to $t$. Since $X$ has
at most $|A| \times | \mo t \mc |$
elements (every element of $X$ can be written in the form $a t^i$ for some
$a \in A$ and $0 \le i < |\mo t \mc |$), this homomorphism is an isomorphism.

\smallskip
(b) The defining relations of $Y$ tell us that every element of $Y$
can be written as $at^i$, with $a\in A$
and $i \in \{0, 1, \dots, 2^q -1\}$.
To prove uniqueness, it is enough to show that
$|Y| = 2^q \cdot |A|$.  Note that $Y$ is the quotient 
of $X = A \sdp \mo t \mc$ by the central cyclic 
subgroup $C = \mo c t^{-2^q} \mc$.  (This subgroup
is central in $X$ because $c^{2^s} = 1$ in $A$.) Since 
$c$ has order $2^{-q} |\mo t \mc|$ in $A$ and $t^{2^q}$ 
has order $2^{-q} |\mo t \mc|$ in $\mo t \mc$, 
we have $|C| = 2^{-q} |\mo t \mc|$
\[ |Y| = \frac{|X|}{|C|} = \frac{|A| \cdot |\mo t \mc|}{| C |} = 
2^q |A| \, , \]
as desired.

\smallskip
(c) Every element of $S$ can be uniquely written
in the form $a t^i$, for some $a \in A$ and $0 \le i < 2^q$. Thus
the natural surjective homomorphism $Z \lra S \simeq \mo A, t \mc$
is an isomorphism.
\end{proof}

We are now ready to prove the main result of this section. We will
continue to use the notations of Lemma~\ref{lem.generators}.

\begin{reduction} \label{red5.2}
In the proof of the implication $(d) \Longrightarrow (a)$
of Theorem~\ref{thm2} we may assume without loss of generality that

\smallskip
(1) $e_1 = \dots = e_r$ and

\smallskip
(2) $S$ is a semidirect product of $A$ and $\mo t \mc$.
\end{reduction}

\begin{proof}
We will use the following two simple ``moves" to
go from an arbitrary Iwasawa group to one satisfying (1) and (2):

\smallskip
(i) If $H$ is a subgroup of $G$ and $G$ admits
a non-hyperbolic trace form then so does $H$.

\smallskip
(ii) Suppose $T$ is a $2$-group and $N$ be a normal subgroup of $T$
contained in $T^2 = \Fr(T)$.
If $T$ admits a non-hyperbolic trace form then so does $T/N$.

\smallskip
\noindent
(ii) is immediate from Corollary~\ref{cor.serre}.
To prove (i), note that if the trace form of a $G$-Galois extension $L/K$
is not hyperbolic then neither is the trace form of $L/L^H$;
see, e.g.,~\cite[Lemma 2.1(c)]{kr}

\smallskip
(1) Let $e = \max \{ e_1, \dots, e_r \}$ and embed $A$ in the abelian
group \[ B = \mo b_1 \mc \times \dots \times \mo b_r \mc \simeq
\bbZ/2^e \times \dots \times \bbZ/2^e \, , \]
where each $b_i$ has order $2^e$ and $a_i = b_i^{2^{e-e_i}}$ for all
$i=1,2,\dots,r$. Let
\[ S_1 = B * \mo t \mc / \mo t^{2^q} = a_0 \, , \;
t b t^{-1} = b^{1 + 2^s} \, | \, b \in B \mc \, . \]
Then there is a natural homomorphism $S \simeq Z \lra S_1$,
which sends $t$ to $t$ and $a$ to $a$ for every $a \in A \subset B$.
By Lemma~\ref{lem.generators}(b),
this homomorphism is injective. Thus by (i) we may
replace $S$ by $S_1$. This completes the proof of (1).

\smallskip
 From now on, we will assume that $e_1 = \dots = e_r = e$.

\smallskip
(2) Let $X$ and $Z$ be as in Lemma~\ref{lem.generators}. Consider the natural
homomorphism $f \colon X \lra Z \simeq S$ which sends $t$ to $t$ and
$a$ to $a$ for every $a \in A$. By Lemma~\ref{lem.generators}(a)
$X \simeq A \sdp \mo t \mc$.
It now suffices to show that $\Ker(f) \subset \Fr(X) = X^2$;
part (2) will then follow from (ii), with $T = X$. For notational
convenience, we will denote the image $t$ in $S$ by $\overline{t}$.

Suppose $a t^i \in \Ker(f)$ for some $a \in A$
and $0 \leq i < | \mo t \mc|$;
in other words, $a \overline{t} \, ^i = 1$ in $S$.
Then, since the order of $\overline{t} \, A$ 
in $S/A$ is $2^q$, we conclude that $i$ is 
a multiple of $2^q$. In particular, since $S$ is not abelian,
we have $q \geq 1$ and thus $t^i \in X^2$. 
It remains to show that $a \in X^2$.
Indeed, since $a = \overline{t} \, ^{-i}$ in $S$,
$a$ and $\overline{t}$ commute in $S$, i.e.,
$a^{2^s} = 1$ in $S$. Since we are assuming that
$A \simeq (\bbZ/ 2^e \bbZ)^r$ and $s < e$,
cf. part (1) and~\eqref{e.s}, we conclude that
$a \in A^2$ in $A$, and consequently $a \in X^2$ in $X$,
as claimed.
\end{proof}

\section{Conclusion of the proof of Theorem~\ref{thm2}
(d) $\Longrightarrow $ (a)}
\label{sect.pf3}

In view of Reduction~\ref{red5.2}, it remains to prove the following

\begin{prop} \label{prop6.1}
Let $S = A \sdp \mo t \mc$, where $\mo t \mc$ is a finite 
cyclic $2$-group, acting on $A=(\bbZ/2^{e}\bbZ)^r$
by $tat^{-1}=a^{1+2^{s}}$, and $2 \leq s < e$.
Then there exists a $S$-Galois extension $E/F$ such 
that $F$ contains a primitive root of unity $\zeta_{2^s}$ 
of degree $2^s$ and the trace form $q_{E/F}$ is non-hyperbolic.
\end{prop}

Our proof of Proposition \ref{prop6.1} below relies on valuation theory;
our primary background references are \cite{ek},~\cite{rib} 
and \cite{wadsworth}.  We shall denote the finite field 
of order $q$ by $\bbF_q$.

\begin{lem} \label{lem.valuation}
For every integer $s \geq 2$,
there exists a field $F$ with a 2-henselian valuation $v$
with value group $\Gamma_v$, and residue field $\mathcal{K}$,
such that

\smallskip
(i) $\Char \, \mathcal{K} \neq 2$,

\smallskip
(ii) $F$ contains a primitive root of unity $\zeta_{2^s}$
of degree $2^s$ but does not contain the primitive root
of unity $\zeta_{2^{s+1}}$ of degree $2^{s+1}$,

\smallskip
(iii) $\dim_{\bbF_2}\Gamma_{v}/2\Gamma_{v} \geq r$.

\smallskip
(iv) $\mathcal{K}(2)=\mathcal{K}(\zeta_{2^{\infty}})$, 
where $\mathcal{K}(\zeta_{2^{\infty}})$ is the extension
of $\mathcal{K}$ obtained by adjoining all $2^n$th roots 
of unity to $\mathcal{K}$, for $n = 1,2,\dots$ and
$\mathcal{K}(2)$ is the maximal 
$2$-extension of $\mathcal{K}$ in some algebraic
closure of $\mathcal{K}$.

\smallskip
\noindent
Moreover, we can choose $F$ so that $\Char(F) = 0$.
\end{lem}

\begin{proof} We shall give two constructions: a simple one in prime
characteristic and a slightly more complicated one in characteristic zero.

\smallskip
Construction 1: 
Observe that $5^{2^{s-2}}-1$ is divisible by $2^s$ but not by $2^{s+1}$
for any integer $s\geq 2$; see, e.g.,~\cite[5.3.17]{scott}.  
Therefore if $q=5^{2^{s-2}}$ then
$\zeta_{2^s}\in \bbF_q$ but $\zeta_{2^{s+1}} \notin \bbF_q$.
Let $F= \bbF_{q}((X_{1}))((X_{2}))\dots((X_{r}))$
be the field of the iterated power series in variables
$X_{1},\dots,X_{r}$ over $\bbF_{q}$ and $v$ be the natural $2$-henselian
valuation $v: F\lra \bbZ \times \dots \times \bbZ$ ($r$-times),
where $\bbZ\times\dots\times\bbZ$ is lexicographically ordered.
One also has $\mathcal{K}(v)= \bbF_{q}$, so that properties (i)-(iv) hold.

\smallskip
Construction 2: Alternatively consider the field
$$F=\mathbb{Q}_p((x_1))((x_2))\dots((x_r))$$
of characteristic $0$ and the natural $2$-henselian valuation
$$ \text{$v \colon F\lra\bbZ \times \dots \times \bbZ$ ($r$ times).} $$
This valuation composed with the $p$-adic valuation on $\mathbb{Q}_p$ 
(see e.g.,~\cite[p. 63]{rib}) yields a new $2$-henselian valuation $v'
:F\lra\bbZ\times\dots\times\bbZ$ ($(r+1)$-times) with a residue field 
$\mathcal{K}(v')=\bbF_p$.
(The fact that $v'$ is again $2$-henselian 
follows from~\cite[Proposition~10, page~211]{rib}; see 
also~\cite[p. 4]{koe}.) Thus $v'$ satisfies conditions 
(i), (iii) and (iv). 

It remains to show that we can choose the prime $p$ so that 
condition (ii) holds.
We claim that for each $s\in\mathbb{N}$ there is
a prime $p$ such that $\zeta_{2^{s}}\in\mathbb{Q}_{p}$ but
$\zeta_{2^{s+1}}\notin\mathbb{Q}_{p}$.
By Hensel's Lemma it is enough to show that
for each $s\in\mathbb{N}$ there exists a prime
$p$ such that $p -1$ is divisible by $2^{s}$ but
not by $2^{s+1}$.  To construct $p$, note that
by Dirichlet's theorem there exists $n\in\mathbb{N}$ such
that $p=(1+2^{s})+ 2^{s+1}n$ is a prime number; this
prime $p$ has the desired properties.
\end{proof}

For the rest of this section, we shall assume that $F$, $v$,
$\Gamma_v$ and $\mathcal{K}$ are as in Lemma~\ref{lem.valuation},
$\bbZ_2$ is the additive group of $2$-adic integers
and furthermore,

\begin{itemize}

\smallskip
\item
$F(2)$ is the maximal $2$-extension of $F$ in some algebraic closure,

\smallskip
\item $G_{F}(2):=\Gal(F(2)/F)$ is the Galois group of $F(2)/F$,

\smallskip
\item $T_{v}\simeq \bbZ_2 \times \dots \times \bbZ_2$ ($d$-times),
where $d =\dim_{F_{2}}\Gamma_{v}/2\Gamma_{v}$. Here $T_{v}$
denotes the inertia subgroup of $G_F(2)$ associated with $v$,

\smallskip
\item
$w$ is the unique valuation of $F(2)$ which extends $v$ on $F$.
\end{itemize}

\smallskip
\noindent
By a result of Engler and Koenigsmann~\cite[Proposition~1.1b]{ek},
\[ G_F(2) \simeq (T_{v}\times G_{\mathcal{K}(\zeta_{2^{\infty}})}(2)) \sdp
\bbZ_{2} \, , \]
where $\bbZ_{2}=<\sigma>$ and the action of $\sigma$ on $T_{v}$ is
$\sigma^{-1}\tau\sigma=\tau^{2^{s}+1}$ for every $\tau\in T_{v}$.

It is also worthwhile to recall that $T_{v}/T_{v}^{2}$ is the Pontrjagin
dual of $\Gamma_{v}/2\Gamma_{v}$, and this
duality is induced by the Kummer pairing
\[ \mo \; \; ,  \; \; \mc \, \colon
T_{v}/T_{v}^{2}\times\Gamma_{v}/2\Gamma_{v}\lra\{\pm 1\} \, , \]
where $\mo  \, [\theta] \, , \, [f] \, \mc=\theta(\sqrt{f})/\sqrt{f}$
for each $\theta\in T_{v}$ and $f\in F^{*}$.
Here $[\theta]\in T_{v}/T_{v}^{2}$ and
$[f]\in\Gamma_{v}/2\Gamma_{v}$ denote the images in
$\theta$ and $f$ in the factor groups $T_{v}/T_{v}^{2}$ and
$\Gamma_{v}/2\Gamma_{v}$, respectively.

We are now ready to finish the proof of Proposition~\ref{prop6.1}.
Suppose $G_{\mathcal{K}(\zeta_{2^{\infty}})}(2)=\{1\}$, i.e.,
$\mathcal{K}(2)=\mathcal{K}(\zeta_{2^{\infty}})$. Then we have
$$G_{F}(2)\simeq T_{v}\sdp \bbZ_{2}.$$
Since
$d=\dim_{F_{2}}\Gamma_{v}/2\Gamma_{v}\geq r$ we deduce that
\[ T_{v}=
\underbrace{\bbZ_2 \times \dots \times \bbZ_2}_{\text{$r$ times}}
\times S \]
for some suitable subgroup $S$ of $T_{v}$.
Therefore there exists a surjective homomorphism
$\tilde{\varphi}:T_{v}\lra A$ which projects the first factor
on $A$ and is trivial on $S$. Because the action
of $\sigma$ on $T_{v}$ is given by $\sigma^{-1}\tau\sigma=\tau^{1+2^{s}}$
for each $\tau\in T_{v}$, we see that
$\tilde{\varphi}$ extends uniquely to a surjective homomorphism
$$\varphi:G_{F}(2)\lra S\mbox{ such that }\varphi(\sigma)=t^{-1}.$$
Let $R$ be the kernel of $\varphi$ and $E$ the fixed
field of $R$. Then $E/F$ is Galois and
$\Gal(E/F)\simeq S$. From the fact
that $T_{v}\simeq \Hom(\Gamma_{w}/\Gamma_{v},\zeta_{2^{\infty}})$
(see \cite[page 2474]{ek}) and
the fact that the outer factor $\bbZ_{2}$ in the semidirect decomposition of
$G_{F}(2)$ as $T_{v}\sdp \bbZ_{2}$
is $\Gal(F(\zeta_{2^{\infty}})/F)$,
we see that the maximal Galois subextension
$E'/F$ of $E/F$ with a Galois group of
exponent $2$ has the form
$$E' =F(\sqrt{a_1},\dots,\sqrt{a_r},\zeta_{2^{s+1}}),$$
where $a_1,a_2,\dots,a_r\in F^{*}$ such that their
values $v(a_1),\dots,v(a_r)\in\Gamma_{v}$ are linearly independent
in $\Gamma_{v}/2\Gamma_{v}$ over $\bbF_{2}$.

 From \cite[Proposition 4.7]{wadsworth} we see that the Pfister form
$$\ll a_1,\dots,a_r,\zeta_{2^{s}}\gg$$
is non-hyperbolic.  By Corollary~\ref{cor.serre}
the trace form of $E/F$ is Witt equivalent
to a scalar multiple of $\ll a_1,\dots,a_r,\zeta_{2^{s}} \gg$,
which is also non-hyperbolic.
This completes the proof of Proposition~\ref{prop6.1} and thus
of Theorem~\ref{thm2}.
\qed

\begin{remark} \label{rem.char0} Our proof
shows that if the equivalent conditions (a) - (d) of Theorem~\ref{thm2}
hold then the fields $F$ and $K$ in parts (a) and (b)
can be chosen to be of characteristic zero.
\end{remark}

\section{Applications}
\label{sect.appl}

\subsection*{Trace forms over ``small" fields}

\begin{prop} \label{prop2} Let $G$ be a finite group,
$S$ be a Sylow $2$-subgroup of $G$,
$K$ be a field containing a primitive
$4$th root of unity and $L/K$ be a $G$-Galois extension.
Denote the Frattini rank of $S$ by $r$.

\smallskip
(a) If $K$ is a $C_{r-1}$-field then the trace form $q_{L/K}$ is
hyperbolic.

\smallskip
(b) If $\cd_2(K) \leq r-1$ then the trace form $q_{L/K}$ is hyperbolic.

\smallskip
(c) If $K$ is a number field and $r \geq 3$ (i.e., $S$ cannot be
generated by two elements) then the trace form $q_{L/K}$ is hyperbolic.
\end{prop}

Here $\cd_2(K)$ refers to the $2$-cohomological dimension of $K$.
For the definition of cohomological dimension and
of the $C_i$ property for fields, see~\cite[II.4]{serre}.

\begin{proof} By Theorem~\ref{thm1} it is enough to show that
under the assumptions of the corollary every
$r$-fold Pfister form $q$ over $K$ is hyperbolic.

In part (a) $q$ is necessarily isotropic and, hence, hyperbolic;
see, e.g.,~\cite[Corollary 10.1.6]{lam}. In part (b),
by Milnor's conjecture (recently proved by Voevodsky~\cite{voevodsky})
$q$ lies in $I^{r+1}$, where $I$
is the fundamental ideal in the Witt ring $W(K)$ and by the Arason-Pfister
theorem this is only possible if $q$ is
hyperbolic; see~\cite[Corollary 10.3.4]{lam}.

Part (c) is a special case of (b), since a totally imaginary number field
has cohomological dimension 2; see \cite[II.4.4]{serre}. However, a much 
more elementary argument, based on the Hasse-Minkowski principle, 
is available in this case. Indeed, every quadratic form of 
dimension $\ge 5$ over $K$ is isotropic; see~\cite[Corollary 3.5, p. 169]{lam}.
In particular, for $r \ge 3$, every $r$-fold Pfister form is isotropic
and hence hyperbolic over $K$.
\end{proof}

\subsection*{Simple groups}

\begin{prop} \label{prop7.8}  Let $G$ be a finite simple group and
let $S$ be the Sylow $2$-subgroup of $G$. Then
the following are equivalent.

\smallskip
(a) $S$ is abelian, and

\smallskip
(b) There exists a $G$-Galois field extension $L/K$ such that
$K$ contains a primitive $4$th root of unity and
the trace form $q_{L/K}$ is not hyperbolic.
\end{prop}

\begin{proof}
By Theorem~\ref{thm2} it is sufficient to prove that $S$ cannot
be a non-abelian Iwasawa group. Equivalently
(via Iwasawa's theorem~\cite{iwasawa}) $S$ cannot be a non-abelian
modular non-Hamiltonian $2$-group. The last assertion
is an immediate consequence of~\cite[Proposition~4.2]{ward}.
(It can also be deduced from~\cite[page~197, Exercise~1]{schmidt}.)
\end{proof}

For the sake of completeness we remark if a finite simple group
$G$ has an abelian $2$-Sylow subgroup $S$ then $S$ is necessarily
elementary abelian (see~\cite[Theorem 4.2.3]{feit}); moreover,
Walter~\cite{walter2} classified all finite simple groups $G$ with
this property.

\subsection*{The extension problem}

Let $G$ be a finite group and $N$ be a normal subgroup of $G$ and
$K \subset L$ be a $G/N$-Galois field extension. Recall that the
{\em extension problem} for this data is the question of existence
of a tower $K \subset L \subset M$, such that $M/K$ is a $G$-Galois
field extension,
and the natural quotient map $\Gal(M/K) \lra \Gal(L/K)$ coincides
with $G \lra G/N$.

Now assume that $G$ is a nonabelian $2$-group of Frattini
rank $r$, $N = \Fr(G) = G^2$, and
$L = K(\sqrt{a_{1}},\dots,\sqrt{a_{r}})$ is a multiquadratic extension of
$K$ of degree $2^r$ such that $\Gal(L/K) \cong G/\Fr(G) = (\bbZ/ 2 \bbZ)^r$.
Assume also that $K$ contains a primitive $e$th root of unity,
where
\[ e = \min \{ \exp(H) \, | \, \text{$H$ is
a non-abelian subgroup of $G$} \} \, . \]

\begin{prop} \label{prop7.1} 
If the extension problem for $G$, $N$, and $L/K$ defined above
has a solution, then the $r$-fold Pfister form
$\ll a_{1},\dots,a_{r}\gg$ is a hyperbolic over $K$.
\end{prop}

\begin{proof} Suppose
$L/K$ is the required $G$-Galois field extension.
Then from Theorem~\ref{thm0} we see that the trace
form $q_{L/K}$ is hyperbolic. But from
Corollary~\ref{cor.serre}(a) we see that $q_{L/K}$
is Witt equivalent to a scalar multiple
of $\ll a_{1},\dots,a_{r}\gg$. Hence
$\ll a_{1},\dots,a_{r}\gg$ is hyperbolic as required.
\end{proof}

\section{Which quadratic forms are trace forms?}

We now return to the question we posed at the beginning
of the Introduction.  Let $G$ be a finite group and $K$
be a field containing $\sqrt{-1}$.  Which quadratic forms $q$
over $K$ can occur as trace forms of $G$-Galois field
extension $L/K$?
In view of Theorem~\ref{thm2} we may assume that
the Sylow $2$-subgroup $S$ of $G$ is an Iwasawa $2$-group;
otherwise every trace form will be hyperbolic. By Theorem~\ref{thm1}
\[ q \simeq |S| \otimes \text{($r$-fold Pfister form)} \]
but, in general, we do not know which $r$-fold Pfister
forms can occur, even if $G = S$ is a $2$-group.
In this section we will describe the trace
forms for one particular family of groups.

Recall that the modular group $M(2^n)$ of order $2^n$ is defined as
$$M(2^{n})= \mo \, \sigma,\tau \, | \, \sigma^{2^{n-1}}=1=\tau^{2},
\tau\sigma\tau=\sigma^{1 + 2^{n-2}} \, \mc \, .$$
In the sequel 
\begin{equation} \label{e.n>= 4}
\text{we will always assume that $n \geq 4$.} 
\end{equation}
It is easy to see that $M(2^n)$ is an Iwasawa group
of order $2^n$, exponent $2^{n-1}$ 
and strength $n-2$.
Setting $A = \mo \sigma \mc$, we see that $(A, \tau)$ is
an Iwasawa structure on $M(2^n)$ of level $n-2$.
Note also that the Frattini subgroup of $M(2^n)$
is $\Fr(M(2^n)) = \mo \sigma^2 \mc$.

\smallskip
For future reference we record the following elementary observation.
As usual, we shall denote the class of $a \in K^*$ 
in $K^{\ast}/(K^{\ast})^2$ by $[a]$.

\begin{remark} \label{rem1.modular}
Let $K$ be a field containing a primitive 
$4$th root of unity $\zeta_4$. Then $2 \zeta_4 = (1 + \zeta_4)^2$ 
and thus 
\begin{equation} \label{e1.modular}
\text{$[2] = [\zeta_4]$ in $K^{\ast}/(K^{\ast})^2$.} 
\end{equation}
In particular,

\smallskip
(i) if $K$ contains a primitive $8$th root of unity then $2$ is
a square in $K$ and

\smallskip
(ii) if $K$ contains a primitive root of unity $\zeta_{2^{n-2}}$
then $2^n$ is a square in $K$.

\smallskip
Indeed, (i) is immediate from~\eqref{e1.modular}.
To prove (ii), consider two cases: $n = 4$ and $n \geq 5$; 
see~\eqref{e.n>= 4}. If $n = 4$ then $2^4 = 4^2$ is
certainly a square.  For $n \ge 5$ (cf.~\eqref{e.n>= 4}), 
(ii) follows from (i).
\end{remark}

We now proceed with the main result of this section. As usual,
$\zeta_i$ will denote a primitive $i$th root of unity. 

\begin{prop} \label{prop.modular}
Let $n \geq 4$ be an integer, $K$ be a field such that 
$\zeta_{2^{n-2}} \in K$ but $\zeta_{2^{n-1}} \not \in K$ and
$q$ be a non-degenerate $2^n$-dimensional quadratic form over $K$.
Then the following are equivalent:

\smallskip
(a) $q$ is Witt equivalent to the trace form 
of some $M(2^n)$-Galois field extension $L/K$. 

\smallskip
(b) $q$ is Witt equivalent to
$\ll \zeta_{2^{n-2}}, a \gg$ for some $a \in K^*$, where 
$[a] \neq [1]$, $[\zeta_{2^{n-2}}]$ in $K^*/(K^*)^2$.
\end{prop}

Our assumption that $\zeta_{2^{n-1}} \not \in K$ is harmless,
since otherwise Theorem~\ref{thm0} tells us that
the trace form of every $M(2^n)$-Galois extension is hyperbolic.
On the other hand, the assumption that $\zeta_{2^{n-2}} \in K$ 
is essential. 

\begin{proof} 
Set $K' = K(\zeta_{2^{n-1}})$, where
$\zeta_{2^{n-1}}$ is a primitive root of unity
of degree $2^{n-1}$. By our assumption on $K$, $[K':K] = 2$.

\smallskip
(b) $\Longrightarrow$ (a): Suppose
$q \simeq \ll \zeta_{2^{n-2}}, a \gg$, where 
$a \neq [1]$, $[\zeta_{2^{n-2}}]$ in $K^*/(K^*)^2$.
We will construct an $M(2^n)$-Galois
extension $L/K$ whose trace form is Witt equivalent to $q$ by
modifying~\cite[Example 6.1]{kr}, due to Serre. 

Let $L = K'(\sqrt[2^{n-1} \; \, \,]{a})$. By our assumption on $[a]$,
$a$ is not a square in $K'$. Thus $[L : K'] = 2^{n-1}$
(see, e.g.,~\cite[Theorem VIII.9.16]{lang}) and consequently,
$[L:K] = 2^n$. Now the computations in~\cite[Example 6.1]{kr}
show that $L/K$ is an $M(2^n)$-Galois extension whose trace form
$q_{L/K}$ is Witt equivalent to $\mo 2^n \mc \otimes
\ll \zeta_{2^{n-2}}, a \gg$. Finally by 
Remark~\ref{rem1.modular}(ii), $2^n$ is a square in $K$ 
and thus the factor of $\mo 2^n \mc$ can be removed.
In other words, $q$ is Witt equivalent to
$\ll \zeta_{2^{n-2}}, a \gg$, as claimed.

(a) $\Longrightarrow$ (b): Assume that $q = q_{L/K}$ for some $M(2^n)$-Galois 
extension $L/K$. Then $q \otimes_K K'$ is the trace form
of the $M(2^n)$-Galois $K'$-algebra $L \otimes_K K'$.
By Theorem~\ref{thm0}, we know that $q \otimes_K K'$
is hyperbolic. (Recall that Theorem~\ref{thm0} applies
to Galois algebras as well as field extensions; see the
first remark after the statement of Theorem~\ref{thm1}
in Section~\ref{sect.intro}.) On the other hand, combining
Theorem~\ref{thm1} and Remark~\ref{rem1.modular}, we see that
$q$ is Witt equivalent to a $2$-fold Pfister form. 
The basic theory of Pfister forms (see, e.g.,~\cite[p. 465]{arason})
now tells us that $q$ is Witt equivalent to
$\ll \zeta_{2^{n-2}}, a \gg$ for some $a \in K^{\ast}$.

It remains to show that $a$ can always be chosen so that
$[a] \neq [1]$, $[\zeta_{2^{n-2}}]$ in $K^*/(K^*)^2$. 
Note that if $[a] = [1]$ or $[\zeta_{2^{n-2}}]$
then $\ll \zeta_{2^{n-2}}, a \gg$ is a hyperbolic trace form.
Thus in order to finish the proof of the proposition, it suffices to 
establish assertions (i) and (ii) below. 
Recall that a field $K$ containing a primitive $4$th root of unity $\zeta_4$ 
is called {\em rigid} if and only if for every $k \not \in (K^*)^2$, the form
$\mo 1, k \mc$ represents only the classes $[1]$ and $[k]$ in $K^*/(K^*)^2$;
cf.~\cite[Section 3]{ware}.

\smallskip
(i) If $K$ is rigid then no $M(2^n)$-Galois field extension
$L/K$ has a hyperbolic trace form. 

\smallskip
(ii) If $K$ is not rigid then $\ll \zeta_{2^{n-2}}, b \gg$ 
is hyperbolic for some $b \in K^*$ such that 
$[b] \neq [1]$, $[\zeta_{2^{n-2}}]$ in $K^*/(K^*)^2$. 

\smallskip
\noindent
In other words, if $K$ is rigid then the case
where $[a] = [1]$ or $[\zeta_{2^{n-2}}]$ can never occur. If
$K$ is not rigid then, after possibly replacing $a$ by $b$, we 
can always assume that $[a] \neq [1]$, $[\zeta_{2^{n-2}}]$ in $K^*/(K^*)^2$. 

\smallskip
To prove (i), note that
if $L/K$ is an $M(2^n)$-Galois extension then
$L^{\Fr(M(2^n))}$ is a $\bbZ/2 \times \bbZ/2$-Galois
extension of $K$. Hence, $L^{\Fr(M(2^n))}$ has the form
$K(\sqrt{a}, \sqrt{b})$ for some $a, b \in K^*$,
where $a$ and $b$ are $\bbF_2$-linearly
independent in $K^{\ast}/(K^{\ast})^2$. By
Corollary~\ref{cor.serre}(a),
\[ q \simeq \mo |\Fr(M(2^n))| \mc \otimes q_{K(\sqrt{a}, \sqrt{b})/K} \, . \]
Here $|\Fr(M(2^n))| = 2^{n-2}$ because 
$\Fr(M(2^n))$ is the cyclic subgroup of $M(2^n)$ generated by
$\sigma^2$. Combining this with formula~\eqref{e.elem-abelian} for
$q_{K(\sqrt{a}, \sqrt{b})/K}$, we obtain 
\[ q \simeq \mo 2^n \mc \otimes \ll a, b \gg \simeq \ll a, b \gg \, ,\] 
where the factor of $\mo 2^n \mc$ can be removed in view of 
Remark~\ref{rem1.modular}(ii). Over a rigid field such a form cannot 
be isotropic, since otherwise $\mo 1, a \mc$ would take on the same value 
as $\mo b \mc \otimes \mo 1, a \mc$, thus making $[a]$ and $[b]$
linearly dependent over $\bbF_2$. This proves (i).

To prove (ii), we appeal to~\cite[Theorem 2.16(2)]{ware}, which 
tells us that over a non-rigid field $K$ the
form $\mo 1, \zeta_{2^{n-2}} \mc$ assumes a value $b$ such that
$[b] \neq [1]$, $[\zeta_{2^{n-2}}]$ in $K^*/(K^*)^2$. Then
$\ll \zeta_{2^{n-2}}, b \gg$ is hyperbolic, as claimed. 
\end{proof}

\begin{remark} \label{rem2.modular}
Suppose $n = 4$. Then by Remark~\ref{rem1.modular},
we can replace the form $\ll \zeta_4, a \gg$ in the statement
of Proposition~\ref{prop.modular} by $\ll 2, a \gg$. This way we
recover~\cite[Corollary 6(b)]{dek} for $G = M(16)$.
\end{remark}

\bibliographystyle{amsplain}

\def\cprime{$'$}
\providecommand{\bysame}{\leavevmode\hbox to3em{\hrulefill}\thinspace}

\end{document}